\documentclass[11pt,leqno]{article}

\usepackage{amssymb,amsmath,amsthm}

\def\F {\bb F}
\newcommand{\n}{\noindent}
\newcommand{\bb}[1]{\mathbb{#1}}
\newcommand{\cl}[1]{\mathcal{#1}}
\newcommand{\vp}{\varepsilon}
\newcommand{\ovl}{\overline}

\theoremstyle{plain}
\newtheorem{thm}{Theorem}
\newtheorem{pro}[thm]{Proposition}
\newtheorem{lem}[thm]{Lemma}
\newtheorem{cor}[thm]{Corollary}

\theoremstyle{remark}
\newtheorem*{rk}{Remark}
\newtheorem{rem}[thm]{Remark}

\theoremstyle{definition}
\newtheorem{defn}[thm]{Definition}

\begin{document}

\title{Simultaneous similarity, bounded generation and amenability}

\author{by\\ Gilles Pisier\footnote{Partially supported by  N.S.F. grant 
No.~0200690}\\
Texas A\&M University\\
College Station, TX 77843, U. S. A.\\
and\\
Universit\'e Paris VI\\
Equipe d'Analyse, Case 186, 75252\\
Paris Cedex 05, France}

\date{}
\maketitle
\begin{abstract}
 We prove that a discrete group $G$  is amenable iff it is strongly unitarizable in the following sense:  every unitarizable representation $\pi$ on $G$ can be unitarized by an invertible chosen in the von Neumann algebra generated by the range of $\pi$. Analogously a $C^*$-algebra $A$ is nuclear iff any bounded homomorphism
$u:\ A\to B(H)$  is strongly similar to a $*$-homomorphism  in the  sense that  there is an invertible operator $\xi$ in the
von Neumann algebra generated by the range of $u$ such that
$a\to \xi u(a) \xi^{-1}$ is a $*$-homomorphism.
An analogous characterization holds in terms of derivations.  
We apply this to answer several questions left open in our previous work concerning the length $L(A\otimes _{\max} B)$  of  the maximal tensor product $A\otimes_{\max} B$ of two unital $C^*$-algebras, when we consider its  generation
by the subalgebras $A\otimes 1$ and $1\otimes B$. We show that if
 $L(A\otimes _{\max} B)<\infty$ either for $B=B(\ell_2)$ or when $B$ is the $C^*$-algebra (either full or reduced) of a non Abelian free group, then $A$ must be nuclear. We also show that $L(A\otimes _{\max} B)\le d$ iff the canonical
 quotient map from the unital free product $A{\ast} B$ onto $A\otimes_{\max} B$ 
remains a complete
 quotient map when restricted to the closed span of the words of length $\le d$.
\end{abstract}
 \section*{Introduction}
 
In 1950, J.\ Dixmier and M. Day proved that any amenable group $G$ is unitarizable, i.e.\ 
any uniformly bounded representation $\pi\colon \ G\to B(H)$ is similar to a 
unitary representation. More precisely there is an invertible operator 
$\xi\colon \ H\to H$ such that $\xi\pi(\cdot)\xi^{-1}$ is a unitary 
representation of $G$. The    proof uses a simple averaging argument from 
which it can be seen  that $\xi$ can be chosen with the additional property that 
$\xi$ commutes with any unitary $U$ commuting with the range of $\pi$.
 Equivalently (see Remark \ref {rem11} below),  $\xi$ can be chosen in the von Neumann algebra generated by $\pi(G)$.  (See 
\cite{NW} for more on this). For convenience, let us say that $\pi$ (resp.\ $G$) 
is strongly unitarizable if it has this additional property (resp.\ if every 
uniformly bounded $\pi$ on $G$ is strongly unitarizable).

It is still an open problem whether ``unitarizable'' implies ``amenable'' (see \cite{P7}).
However, we will show that $G$ is amenable iff it is strongly unitarizable. 
Moreover, we will show an analogous result for $C^*$-algebras, as follows.

\begin{thm}\label{thm1}
The following properties of a C$^*$-algebra $A$ are equivalent.
\begin{itemize}
\item[(i)] $A$ is nuclear.
\item[(ii)] For any c.b.\ homomorphism $u\colon \ A\to B(H)$ there is an 
invertible operator $\xi$ on $H$ belonging to the von Neumann algebra generated 
by $u(A)$ such that $a\to \xi u(a) \xi^{-1}$ is a $*$-homomorphism.
\item[(iii)] For any $C^*$-algebra $B$, the pair $(A,B)$ has the following 
simultaneous similarity property:\ for any pair $u\colon \ A\to B(H)$, $v\colon 
\ B\to B(H)$ of c.b.\ homomorphisms with commuting ranges there is an invertible 
$\xi$ on $H$ such that both $\xi u(\cdot)\xi^{-1}$ and $\xi v(\cdot)\xi^{-1}$ 
are $*$-homomorphisms.
\item[(iv)] Same as (iii) but with $v$ assumed to be itself a $*$-homomorphism.

\end{itemize}
\end{thm}

\begin{rem}\label{rem1}
It is possible that (iii) or (iv) for a fixed given $B$ implies that 
$A\otimes_{\min} B = A\otimes_{\max} B$ but this is not clear (at the time of 
this writing).
\end{rem}

\begin{cor}\label{cor2}
If a discrete group $G$ is strongly unitarizable then $G$ is amenable.
\end{cor}

\begin{proof}
Let $A= C^*(G)$. Any bounded homomorphism $u\colon \ A\to B(H)$ restricts to a 
uniformly bounded representation $\pi$ on $G$. Note that $\pi(G)$ and $u(A)$ generate the same 
von~Neumann algebra $M$. Thus if $G$ is strongly amenable, $A$ satisfies (ii) in 
Theorem \ref{thm1}, hence is nuclear and, as is well known, this implies $G$ 
amenable in the discrete case (see \cite{L}).
\end{proof}

 Actually, we obtain a stronger statement:
\begin{cor}\label{cor3}
 If every {\em unitarizable\/} representation  $\pi$ on 
a discrete group $G$ is strongly unitarizable then $G$ is amenable.
 \end{cor}
 \begin{proof}
  Indeed, in Theorem 
\ref{thm1}, $u$ is assumed c.b. on $A=C^*(G)$, so the corresponding $\pi$ is 
unitarizable.
\end{proof}

 \begin{rem}\label{rem11} 
 Assume $G$ amenable with invariant mean
 $\phi$.  Consider  a uniformly bounded representation  $\pi$ on $G$,
 then, the  proof  of  the Day-Dixmier  theorem  is as follows:  Essentially we can define $\xi$ by the (non rigorous) formula
$$\xi=(\int \pi(g)^*\pi(g) \phi(dg))^ {1/2}.$$
 But it is obvious how to make this rigorous:  for any $h$ in $H$
 we define $x_h(g)= \langle \pi(g) h, \pi(g)h \rangle $ (note
 $x_h \in L_\infty(G)$) and then define $\xi$
 by setting $\langle\xi^2 h, h \rangle = \phi (x_h)$. Clearly
 (by the invariance of $\phi$)   $\xi$ unitarizes $\pi$, and the above formula  makes it clear
 that  $\xi$ is in the von Neumann algebra generated by the range of $\pi$.
 \end{rem}
 
\begin{rem}\label{rem2} Note that, by \cite{Co}, $C^*(G)$ is nuclear for any separable, 
   connected locally compact group
 $G$, hence every continuous unitarizable representation on $G$ is strongly unitarizable;
 therefore  we definitely must restrict the preceding 
Corollary \ref{cor3} to the discrete case.
\end{rem}
 
\begin{rem}\label{rem4}
The following elementary fact will be used repeatedly:\ let $U$ be a unitary 
operator on $H$ and let $\xi\ge 0$ be an invertible on $H$ such that $\xi U 
\xi^{-1}$ is still unitary. Then $\xi U \xi^{-1} = U$. Indeed, we have 
$(\xi U \xi^{-1})^* (\xi U \xi^{-1})=I$ hence $U^*\xi^2U = \xi^2$. 
Equivalently $\xi^2$ commutes with $U$ and hence $\xi = \sqrt{\xi^2}$ also 
commutes with $U$.
\end{rem}

 The above results are proved in the first part of the paper. The second part is devoted to the length
 of a pair of (unital) $C^*$-algebras $A,B$, introduced in \cite{P4} and denoted below by $L(A\otimes _{\max} B)$. 
  Let $W_{\le d}$ be the closed span of the words of length $\le d$ in the unital free product
$A{\ast} B$.
 We will prove that $L(A\otimes _{\max} B)\le d$ iff the restriction to $W_{\le d}$ of the canonical
 quotient map from  $A\ast B$ onto $A\otimes_{\max} B$   is a complete
 quotient map (i.e. it yields a complete isomorphism after passing to the quotient by the kernel).
  This gives a more satisfactory reformulation of the definition in \cite{P4}. 
To establish this, we need to prove that $W_{\le d}$ decomposes naturally (completely isomorphically)
 into a direct sum of Haagerup tensor products  of copies
 of $A$ and $B$ of order $0\le j\le d$ (see Lemma \ref{lem33}).
 The latter result seems to be of independent interest.

\section*{Notation and Background}

 While the first part uses mostly basic $C^*$-algebra theory and c.b. maps
 (for which we refer to \cite{T, Pa}), the second one requires more background from 
 operator space theory, e.g. the Haagerup tensor product, and its connection with free
 products of operator algebras for which we refer
 the reader to \cite {P1} (see also \cite{BLM, ER}).
 
Recall that a linear map $v\colon \ Y\to X$ between operator spaces is called
completely bounded (c.b. in short) if the maps
$v_n=Id\otimes v\colon\ M_n(Y)\to M_n(X)$ are uniformly bounded,
and we set
$$\|v\|_{cb}=\sup_n \|v_n\|.$$

Equivalently, if we denote $K(Y)=K\otimes_{\min} Y$, we have
$\|v\|_{cb}=\|Id\otimes v\colon\ K(Y)\to K(X)\|.$

A c.b. map $v$ is called a
``complete surjection" (or a ``complete quotient map") if there is a constant $c$ 
such that for any
$n\ge 1$, and  any $x$ in $M_n(X)$ with $\|x\|<1$ there is $y$ in $M_n(Y)$ such that
$[x_{ij}]  = [q(y_{ij})]$ with $\|y\|<c$. 

When this holds with $c=1$ we say that $v$ is a 
complete metric surjection.  Note that,  when $n=1$,  any surjection  satisfies this
for some $c>0$. When this holds for $c=1$ (and  for $n=1$), we say 
that $v$ is a metric surjection.
Equivalently, $v$ is a 
complete (resp. metric) surjection   iff $Id\otimes v \colon \ K(Y)\to K(X)$ is
a (resp. metric) surjection.

Let $A,B$ be unital $C^*$-algebras and let $u\colon \ A\to B(H)$ and $v\colon \ B\to
B(H)$ be linear maps. We denote by $u\cdot v\colon \ A\otimes B\to B(H)$ the linear map
defined on the {\sl algebraic\/} tensor product $A\otimes B$ by $u\cdot v(a\otimes b) =
u(a)v(b)$.

We will say (for short) that $u\cdot v$ is c.b.\ on $A \otimes_{\rm min} B$ (resp.\ $A\otimes_{\rm max} B$) if $u\cdot v$ extends to a c.b.\ map on $A \otimes_{\rm min} B$ (resp. $A\otimes_{\rm max} B$). We will use a
similarly shortened terminology for ordinary boundedness instead of the complete one.  

Now assume that $u,v$ are unital homomorphisms with commuting ranges. Then $u\cdot v$ is
a homomorphism on the incomplete algebra $A\otimes B$. By [H], $u\cdot v$ is c.b.\ on $A
\otimes_{\rm max} B$ iff there is an invertible $\xi$ in $B(H)$ such that $\xi
u(\cdot)\xi^{-1}$ and $\xi v(\cdot)\xi^{-1}$ are both $*$-homomorphisms. More
precisely, we have
\begin{equation}\label{eqH}
\|u\cdot v\colon \ A\otimes_{\rm max} B \to B(H)\|_{cb} 
= \inf\{\|\xi\| \|\xi^{-1}\|\}  
\end{equation}
where the infimum runs over all $\xi$ satisfying this.

Now   assume that $v$ is a unital $*$-homomorphism. In that case, we have
\begin{equation}\label{eqHH}
\|u\cdot v\colon \ A\otimes_{\rm max} B\to B(H)\|_{cb} =
 \inf\{\|\xi\| \|\xi^{-1}\|\}\end{equation}

where the infimum runs over all $\xi$ in
$v(B)'$ such that $a\to \xi u(a)\xi^{-1}$ is a $*$-homomorphism. Indeed, this is an
immediate consequence of Remark \ref{rem4} (since $b\to\xi v(b)\xi^{-1}$ is a
$*$-homomorphism iff it maps unitaries to unitaries).

Let $r\colon \ A\to B(H)$ and $\sigma\colon \ B\to B({\cal H})$ be unital
$*$-homomorphisms, and let $\pi\colon \ A\otimes_{\rm min} B\to B(H\otimes {\cal H})$
be their tensor product, i.e.\ $\pi(a\otimes b) = r(a) \otimes \sigma(b)$.

By an $r$-derivation $d\colon \ A\to B(H)$ we will mean 
a derivation with respect to $r$ (i.e.\ $d(a_1a_2) = r(a_1) d(a_2) + d(a_1)r(a_2))$.
Let $r_1\colon \  A\to
B(H\otimes {\cal H})$ and $\sigma_1\colon \  B\to
B(H\otimes {\cal H})$ be the representations defined by $r_1(a) = a\otimes
I$ and
$\sigma_1(b) = I\otimes
\sigma(b)$.
Let $\delta \colon \ A\to B(H\otimes {\cal H})$ be an $r_1$-derivation  such that 
$\delta(A) \subset (I\otimes \sigma(B))'$.  It is easy to check that $\delta \cdot
\sigma_1$ is a $\pi$-derivation on the (incomplete) algebra $A\otimes B$. For any $T$
in $B(H\otimes {\cal H})$, we denote
$$\delta_T(a) = (r(a)\otimes 1)T - T(r(a)\otimes 1).$$
By a result due to E.~Christensen (\cite{C}) we have then
\begin{equation}\label{eqC}\|\delta\cdot \sigma_1\colon \ A\otimes_{\rm min}B 
\to B(H\otimes {\cal H})\|_{cb} = 2\inf\{\|T\|\mid T\in \sigma_1(B)',
\ \delta=\delta_T\}.\end{equation}
 Actually in the present special situation, this is
also equal to the c.b.\ norm of $\delta\cdot\sigma_1$ on $A\otimes_{\rm max} B$.
Indeed, let $C_{\rm max}$ be the latter c.b.\ norm. Then, Christensen's result
implies that the above
$2\inf\|T\|$ is $\le C_{\rm max}$, but since $\pi$ is {\it assumed continuous\/}  on
$A\otimes_{\rm min} B$ it follows that the c.b.\ norm of $\delta_T$ on $A\otimes_{\rm
min} B$ is $\le 2\|T\|$.

 \section*{Proof of the main result}

Actually we will prove a more general result than the above Theorem \ref{thm1}. 
Indeed, we will show that it suffices for $A$ to be nuclear that (iii) or (iv) holds
for a ``large enough'' $C^*$-algebra $B$. It may be that any non-nuclear $B$ can be
used but we can't prove this. Instead we introduce the notion
of a  ``liberal'' $C^*$-algebra which is close to
being the ``opposite'' of nuclearity.

To describe this, we need to introduce the following notion.

We denote by $E^n_\lambda$ the operator space that is the linear span of 
$\lambda(g_1),\ldots, \lambda(g_n)$ in  the von Neumann algebra
   generated by the 
 left regular representation of  the free group $\F_n$ 
(recall that $ g_1,\ldots,  g_n$
 denote the free generators of $\F_n$).
\medskip

\n {\bf Note:}\ We could use $R_n\cap C_n$ instead of $E^n_\lambda$
 (see \cite[p.\ 184]{P1}), but it is 
easier to see the connection with the preceding argument using $E^n_\lambda$.

\begin{defn}\label{def5}
We say that $\{E^n_\lambda\}$ factors uniformly through an operator space $B$ 
if, for any $n\ge 1$, there are mappings
\[
v_n\colon \ E^n_\lambda \to B,\quad w_n\colon \ B\to E^n_\lambda
\]
such that $w_nv_n = id$ and $\sup\limits_n \|v_n\|_{cb} \|w_n\|_{cb} <\infty$.  
\end{defn}
 
 \begin{defn}\label{def6} We  say that a $C^*$-algebra  $B$ is   ``liberal" if it admits a representation
 $\sigma\colon\ B\to B(H)$ such that  $\{E^n_\lambda\}$ factors uniformly through the commutant $\sigma(B)'$. 
 
 \end{defn}

\begin{rk}
Examples of liberal $C^*$-algebras
 are $C^*(\F_\infty)$, $C^*_\lambda({\bb F}_\infty)$ or the von~Neumann algebra
generated by $\lambda({\bb F}_\infty)$. This follows
from \cite[Th. 4.1]{HP}. A fortiori
 (since $\F_\infty$ embeds in  $\F_2$), the same is true for $\F_n$  
for any $n\ge2$.
 By \cite[p.~205]{A} and \cite{W}, $B(\ell_2)$ or the Calkin algebra
$B(\ell_2)/K(\ell_2)$ are liberal.  Clearly, since nuclear passes to quotients, 
any liberal $C^*$-algebra is
non-nuclear by 
\cite[Th. 4.1]{HP}.  
 Apparently, there is no known counterexample to the converse.
\end{rk}
  
\begin{thm}\label{thm1B}
Assume $A\subset B(H)$ and let $B$ be a liberal (unital) $C^*$-algebra. 
The properties in Theorem \ref{thm1} are equivalent to:
\item[(v)] For any $*$-homomorphism $\sigma\colon \ B\to B({\cal H})$ 
there is a constant $C$ such that for any c.b.\ derivation $\delta\colon \ A\to
B(H)\otimes \sigma(B)'$ (relative to the embedding $A\simeq A\otimes I$) for which the
associated $\delta\cdot\sigma_1$ is also c.b.\ on $A\otimes_{\rm min}B$, there is an
operator $T$ in the von~Neumann algebra generated by $A\otimes I$ and $\delta(A)$ such
that
$$\|T\| \le C\|\delta\|_{cb}\quad \hbox{and}\quad \delta(a) = aT-Ta \quad (a\in A).$$
\item[(vi)] There   is a non-decreasing function $F\colon \ 
{\bb R}_+ \to {\bb R}_+$   satisfying the following:\
 For any c.b.\ homomorphism $u\colon \ A\to B(H)$ and any $*$-homomorphism $v\colon \
B\to B(H)$, with commuting ranges such that $u\cdot v\colon \ a\otimes b \to u(a)v(b)$
extends to a c.b.\ homomorphism on $A\otimes_{\rm min} B$, there is an invertible $\xi$
in $v(B)'$ with $\|\xi\| \|\xi^{-1}\| \le F(\|u\|_{cb})$ such that $\xi
u(\cdot)\xi^{-1}$ is a $*$-homomorphism. 
\end{thm}

\begin{proof}[Proof of Theorems \ref{thm1} and \ref{thm1B}]
We will assume $A$ and $u$ unital for simplicity. (iii) $\Rightarrow$ (iv) is trivial.
(iv) 
$\Rightarrow$ (ii) is easy. Indeed, let $M$ denote the von~Neumann algebra 
generated by $u(A)$. Let $v\colon \ M'\to B(H)$ be the inclusion mapping. If 
(iv) holds, we can find $\xi$ such that, for any unitary pair $a,b$ in  $A$, $\xi 
u(a)\xi^{-1}$ and $\xi v(b)\xi^{-1}$ are both unitary. By polar decomposition of 
$\xi$ we may assume $\xi>0$. Then by the preceding remark, $\xi$ must commute 
with $v(M') = M'$, and hence $\xi\in M'' = M$. The implication (i) $\Rightarrow$ 
(iii) follows from the basic properties of the so-called $\delta$-norm, as 
presented in \cite[Th.\ 12.1]{P1} (see also \cite{LM}). Indeed, let $u,v$ be as in (iii).
Clearly, the mapping $v\cdot u\colon \ b\otimes a\to v(b)u(a)$ satisfies 
$\|v\cdot u\colon \ B\otimes_\delta A\to B(H)\|_{cb}\le \|v\|_{cb} 
\|u\|^2_{cb}$. Now, if $A$ is nuclear, $B\otimes_\delta A = B \otimes_{\min} A = 
B \otimes_{\max} A$, and hence $v\cdot u$ defines a c.b.\  homomorphism $\rho$ on 
$B \otimes_{\max} A$. By \cite{H}, there is an invertible $\xi$ such that 
$\xi\rho(\cdot)\xi^{-1}$ is a $*$-homomorphism, from which we conclude that 
(iii) 
holds.

For Theorem \ref{thm1}, it remains only to prove that (ii) implies (i).
We will show (ii)$\Rightarrow$(vi)$\Rightarrow$(v)$\Rightarrow$(i).

Let $u$ be as in (ii). Let $M_u$ be the von~Neumann algebra generated 
by $u(A)$. We first claim that the mapping $\widehat u\colon \ x\otimes y\to u(x)y$ 
extends to a c.b.\ homomorphism from $A\otimes_{\max} M_u'$ to $B(H)$. Indeed, 
since 
$\xi\in M_u$, $\xi\widehat u(\cdot)\xi^{-1}$ is a $*$-homomorphism on $A\otimes
M_u'$. 
By a routine direct sum  argument, (ii) implies  that there is a 
non-decreasing function $F\colon \ 
{\bb R}_+ \to {\bb R}_+$ such that for all $u$ as in (ii), we have
\begin{equation}\label{eq1}
\|\widehat u\|_{cb} \le F(\|u\|_{cb}).
\end{equation} By \eqref{eqHH},
this clearly shows that (ii) implies (vi).

We now show that (vi) implies (v). Assume (vi). Let $\delta,\sigma$ be as in (v).
 Let $u\colon\ A \to M_2(B(H\otimes {\cal H}))$ be  the homomorphism
\[
u\colon \ a\to \begin{pmatrix} a\otimes 1&\delta(a)\\ 0&a\otimes 1\end{pmatrix} 
\in M_2(B(H\otimes {\cal H})),
\]
and define $v \colon\ B \to M_2(B(H\otimes {\cal H}))$ by
\[
u\colon \ b\to \begin{pmatrix}   \sigma_1(b)&0\\ 0& 
\sigma_1(b)\end{pmatrix} 
\in M_2(B(H\otimes {\cal H})),
\]
Note that   $u,v$  have commuting ranges and $\|u\|_{cb}\le
1+\|\delta\|_{cb}$; also $u \cdot v$ is c.b.  on $A\otimes_{\min} B$
because  we assume in (v) that it is so for $\delta\cdot \sigma_1$
(and the representation $\pi=r_1\cdot \sigma_1$ is continuous, hence
c.b.  on $A\otimes_{\min} B$). Then, since we assume (vi), we 
obtain that there is an invertible $\xi$ in $VN(u(A))$ with $\|\xi\| 
\|\xi^{-1}\| \le F(\|u\|_{cb})$ such that $\xi u(\cdot)\xi^{-1}$ is a 
$*$-homomorphism. Reviewing an argument of Paulsen (see either \cite{Pa} or 
\cite[{p. 80}]{P3}) we find that there is an operator $T$ with $\|T\| \le 
2(\|\xi\| \|\xi^{-1}\|)^2 \le 2 F(\|u\|_{cb})$ in the von~Neumann algebra 
generated by $A\otimes 1$ and $\delta(A)$ such that $\delta(a) = [a\otimes 1, 
T]$ for all $a$ in $A$. 
To deduce (v), by homogeneity, we may assume  $\|\delta\|_{cb} =1$
then $\|u\|_{cb} \le 2$ and we find  (v) with  $C= 2F(2)^2$.
 To complete the proof,
it remains to show (v) implies (i), i.e. that (v) implies that $A$ is nuclear.

Let ${\cal H}=\ell_2({\bb F}_\infty)$. Let $W\subset B({\cal H})$ be the von Neumann
algebra generated by the left regular representation 
$\lambda$ on the free group ${\bb F}_\infty$ with $n$ generators, denoted by 
$g_1,\ldots, g_n,\ldots$.    We will first show
that (v) implies (i) in the particular case $B=W'$.
Let $r\colon \ 
A\to B(H)$ be a $*$-homomorphism. By the well known Connes--Choi--Effros results 
(see \cite{P1}), it suffices to show that $r(A)''$ is always injective.
 For 
simplicity, we replace $A$ by $r(A)$. Thus, it suffices to show that $N=A''$ 
is injective or equivalently that $N'$ is injective. By \cite[{Th.2.9}]{P3}, $N'$ is 
injective iff there is a constant $\beta$ such that $\forall n \forall y_i\in N'$ 
there are elements $a_i,b_i\in N'$ with $y_i=a_i+b_i$ such that
\begin{equation}\label{eq3}
\left\|\sum a_ia^*_i\right\|^{1/2} + \left\|\sum b^*_ib_i\right\|^{1/2} 
\le \beta\|(y_i)\|_{R+C}
\end{equation}
where
\[
\|(y_i)\|_{R+C} = \inf\left\{\left\|\sum \alpha_i\alpha^*_i\right\|^{1/2} + 
\left\|\sum \beta^*_i\beta_i\right\|^{1/2}\right\}
\]
where the infimum runs over all the possible decompositions $y_i= 
\alpha_i+\beta_i$ with $\alpha_i,\beta_i$ in $B(H)$.

Note:\ when there is a c.b.\ projection $P\colon \ B(H)\to N'$, then we can take 
$a_i = P\alpha_i$, $b_i=P\beta_i$ and $\beta=\|P\|_{cb}$.

To prove \eqref{eq3}, we first consider an $n$-tuple $(y_i)$ in $N'$ with 
$\|(y_i)\|_{R+C}<1$, so that $y_i=\alpha_i+\beta_i$ with 
\[
\left\|\sum \alpha_i\alpha^*_i\right\|<1 \quad \text{and}\quad \left\|\sum 
\beta^*_i\beta_i\right\| <1.
\]
 We introduce the derivation $\delta\colon \ A\to B(H) \otimes 
W$ defined as follows:
\[\forall a\in A\quad 
\delta(a) = \sum\nolimits^n_1 [a,\alpha_i] \otimes \lambda(g_i).
\]
It is well known (see \cite{HP} or     \cite[p. 185]{P2}) that there is a decomposition 
$\lambda(g_i) = s_i+t_i$ with $s_i,t_i \in B({\cal H})$ satisfying
 \[
\left\|\sum s^*_is_i\right\|^{1/2}\le 1\quad \text{and}\quad \left\|\sum 
t_it^*_i\right\|^{1/2}\le 1.
\]
Thus, since $y_i\in A'$ (and hence $[a,\alpha_i] = -[a,\beta_i]$) we have
\[
\delta(a) = [a\otimes 1, \theta],
\]
where $\theta = \sum^n_1 \alpha_i\otimes s_i - \beta_i\otimes t_i$. Therefore
\[
\|\delta\|_{cb} \le \|\theta\| \le \left\|\sum \alpha_i\alpha^*_i\right\|^{1/2} 
\left\|\sum s^*_is_i\right\|^{1/2} + \left\|\sum \beta^*_i\beta_i\right\|^{1/2} 
\left\|\sum t_i t^*_i\right\|^{1/2} \le 2.
\]
Note that $ \delta\cdot \sigma_1 $ is a finite sum of maps
that are obviously c.b.\ on $A\otimes_{\min} B$.
Since we assume (v),     there is an operator $T$ with $\|T\| \le 
C\|\delta\|_{cb}\le 2C $ in the von~Neumann algebra 
generated by $A\otimes 1$ and $\delta(A)$ such that $\delta(a) = [a\otimes 1, 
T]$ for all $a$ in $A$.
 Note that $T$ belongs to 
$B(H)\ovl\otimes W$. Let $Q\colon \ W\to W$ be the orthogonal projection onto 
the span of $\lambda(g_1),\ldots, \lambda(g_n)$. It is known (see e.g.\ 
\cite[{p. 184}]{P1}) that $\|Q\|_{cb}\le2$, hence if we set
\[
T_1 = (1\otimes Q)(T)\quad \text{we have}\quad \|T_1\| \le 4C 
\]
and moreover since $ \delta(a) = (1\otimes Q) (\delta(a))=(1\otimes Q) [a\otimes 1,
T] =  [a\otimes 1, T_1]$ we have
\[ 
[a\otimes 1, \theta] = [a\otimes 1,
T_1].
\]
We can write $T_1 = \sum z_i\otimes \lambda(g_i)$. We have 
\begin{equation}\label{eq31}
\max\left\{\left\|\sum z_iz^*_i\right\|^{1/2}, \left\|\sum 
z^*_iz_i\right\|^{1/2} \right\} \le \|T_1\| \le 4C
\end{equation}
and since 
$
  \sum[a,\alpha_i] \otimes \lambda(g_i) = \sum[a,z_i]\otimes 
\lambda(g_i)
$ we find
  $\alpha_i-z_i\in A'$.
To conclude, we set
\[
a_i=\alpha_i-z_i,\quad b_i=\beta_i+z_i
\]
we have $a_i\in A'$, $b_i = y_i-a_i\in A'$ and moreover by \eqref{eq31}
\[
\left\|\sum a_ia^*_i\right\|^{1/2} \le \left\|\sum 
\alpha_i\alpha^*_i\right\|^{1/2} + \left\|\sum z_iz^*_i\right\|^{1/2} \le 1+4C
\]
and similarly
\[
\left\|\sum b^*_ib_i\right\|^{1/2} \le 1+4C.
\]
Thus we obtain \eqref{eq3} with $\beta=2(1+4C   )$, which proves that $A$ is nuclear.
This completes the  proof that (v) implies (i) in the case $B=W'$. But if $B$ is
liberal,
  the preceding argument extends:
we replace $W$ by  $\sigma(B)'$
 and $\lambda(g_i)$ by the  elements in  $\sigma(B)'$ corresponding to the basis
 of $E^n_\lambda$. We skip the easy details.
 \end{proof}

  \begin{rk} The preceding proof obviously 
  shows that $A$ is nuclear iff
  
  \begin{itemize}
\item[(vii )] For any $*$-homomorphism $\sigma:\ A \to B(H)$ and any c.b. 
$\sigma$-derivation $\delta
 :\ A \to B(H)$  (i.e. $\delta(ab)=\delta(a)\sigma(b) + \sigma(a)\delta(b)$),
 there is $T$ in the von Neumann algebra generated 
 by the ranges of $\sigma$ and  $\delta$ such that
 $\delta(a)= \sigma(a) T- T\sigma(a)$ for all $a$ in $A$.
    \end{itemize}
   Note that  any nuclear $A$ is amenable (\cite{H2}) and hence has
   a virtual diagonal, i.e. there is a net $t_i
   =\sum_k a_k(i)\otimes b_k(i)$ bounded in $A \widehat \otimes A$ such that for any $a$ in $A$, $a. t_i -t_i .a $ tends to zero in
   $A \widehat \otimes A$. Then, it is easy to see that if $\delta$ is as above,
   and if $T$ is  a weak$*$-cluster point of
   $\sum_k \delta (a_k(i))\sigma(b_k(i))$, we have 
   $\delta(.)= \sigma(.) T- T\sigma(.)$, and $T$ lies
   manifestly in the von Neumann algebra (or even in the weak closure of the
algebra) generated
   by $\sigma(A)\cup \delta(A)$. This remark should be compared with what is known on ``strongly amenable" $C^*$-algebras
   (a smaller class than the nuclear ones), for which we refer the reader
   e.g. to \cite[Section 1.31]{pat}.
    
    \end{rk}
  \begin{rk}
 In the group case the above argument should be compared with  \cite{BF}.
  \end{rk}
\begin{rem}\label{rem8}
An alternate argument (more direct but somewhat less ``constructive'')
for (v)$\Rightarrow $(i) can be 
obtained as in the following sketch. We use the same notation as in the preceding
proof. We argue that
$\delta$ has range into
$B(H) 
\ovl\otimes W$ and note that the latter commutes with $1\otimes W'$. 
Assume $\|\delta\|_{cb}=1$. Then our 
assumption (v) implies that the map $  \delta\cdot \sigma_1$ extends to a 
bounded  map on $A\otimes_{\min} W'$ with norm $\le 2\|T\|\le 2C$. Then
for any 
$x = \sum\limits_{t\in {\bb F}_\infty} x(t)\otimes \rho(t)$ in $A\otimes W'$ we have
\begin{equation}\label{eq4}
\left\|\sum_{t\in {\bb F}_\infty} \delta(x(t)) (1\otimes \rho(t))\right\| \le 
2C\|x\|_{\min}.
\end{equation}
Composing the operator on the left of \eqref{eq4} with $id\otimes \varphi$ where 
$\varphi(T) = \langle T\delta_e,\delta_e\rangle$, we find 
\begin{equation}\label{eq44}
\left\|\sum_i [x(g_i),\alpha_i]\right\| \le C\|x\|_{\min}.
\end{equation}
Note that 
\begin{equation}\label{eq5}
\max\left\{\left\|\sum x(g_i)  x(g_i)^* \right\|^{1/2}, \left\|\sum 
x(g_i)^*x(g_i)\right\|^{1/2} \right\} \le  \|x\|_{\min}  .
\end{equation}
Hence by Cauchy-Schwarz
\begin{equation}\label{eq6}
 \left\|\sum  \alpha_i  x(g_i)  \right\|  
    \le  \|x\|_{\min}  ,
\end{equation}
and
\begin{equation}\label{eq7}
 \left\|\sum     x(g_i)\beta_i  \right\|  
    \le  \|x\|_{\min} .
\end{equation}
Thus \eqref{eq44} implies
 \begin{equation}\label{eq8}
\left\|\sum_i  x(g_i) y_i \right\| \le (C+2)\|x\|_{\min}
\end{equation}
  Letting $x=\sum x_i \otimes \rho(g_i)$, we find that for
all $(x_i)$ in $A$ we have
\[ \left\|\sum_i  x_i  y_i \right\| \le 2(C+2) 
\max\left\{\left\|\sum x^*_ix_i\right\|^{1/2}, \left\|\sum 
x_ix^*_i\right\|^{1/2}\right\} .\]
 Clearly, this remains valid for all $(x_i)$ in $A''=N$ and hence, by \cite[Cor.
5]{P6}  
$N'$ is injective.
\end{rem}

\begin{rem}\label{norm}   Actually, the preceding argument
shows that  $A$ is nuclear if there is a constant
$C$ such that the ordinary norm of
 $\delta\cdot \sigma_1$ on $ A\otimes_{\min} W'$
is  $\le C\|\delta\|_{cb}$.
 \end{rem}

 \section*{Length for a pair of $C^*$-subalgebras}

Let $G_1,G_2$ be two subgroups generating a group $G$.
One says that $G_1,G_2$ generate $G$ with bounded length
(more precisely with length $\le d$) if every element
in $G$ can be written as a product of a bounded
  number of elements either in $G_1$ or in $G_2$
(resp. a product of   at most $d$ such elements).
Equivalently, let $\psi \colon\ G_1 * G_2 \to G$ be the canonical
 homomorphism from the free product onto $G$, then generation
with length $\le d$ is the same as saying that the restriction
of $\psi$ to the subset formed by all the words of length $\le d$
is surjective.

 It turns out
there is a natural analogue of this
in the $C^*$-algebra (or operator algebra) context,
 already considered in \cite{P4}, 
as follows.

Let $A,B$ be $C^*$-subalgebras of a $C^*$-algebra $Z$.
 By convention, we will view
$M_n(A)$ and $M_n(B)$ as subalgebras of $M_n(Z)$ , so that
if $x_1\in M_n(A)$ and $x_2\in M_n(B)$, then
the product $x_1x_2$ belongs to $M_n(Z)$ , and similary
for a product of rectangular matrices.

Now let $d\ge 1$ be an integer. 
We will say that $L(Z;A,B)\le d$ 
or more simply (when there is no ambiguity) that $L(Z)\le d$
if there is a constant $C$ such that for any 
$n$ and any $x$ in  $M_n(Z)$ with $\|x\|_{M_n(Z)} < 1$ 
and for any $\vp>0$ there is an integer $N$ 
for which we can find matrices $x_1,x_2,\ldots, x_d$  and $y_1,y_2,\ldots, y_d$, 
with entries    either 
all in $A$ or all in $B$,
where $x_1,x_2,\ldots, x_d$ are 
of size respectively $n\times N$, $N\times N,\ldots, N\times N$ and $N\times n$, 
and similarly for $y_1,y_2,\ldots, y_d$,
satisfying  $\prod\limits^d_1 \|x_j\|+\prod\limits^d_1 \|y_j\| <C$ 
 and finally such that
\begin{equation}\label{eq100}
\left\|x - \prod^d_1 x_j-\prod^d_1 y_j\right\|_{M_n(Z)} < \vp.
\end{equation}
If this holds but only for $n=1$, then we say that $L_1(Z;A,B)\le d$
or simply that $L_1(Z)\le d$.
Note that the two products are needed because one of them
``starts" in $B$ and the other ``starts" in $A$. So we will make
the convention that $x_1$ is a matrix with entries in $B$
while  $y_1$ is one  with entries in $A$. To
eliminate the $\vp$-error term, we need to use infinite matrices, as follows. 
Let us denote $K (A) = K \otimes_{\min} A$. We may identify $K(A)$ and $K(B)$ 
with subalgebras of $K(Z)$. Then $L(Z) \le d$ (resp.\ 
$L_1(Z)\le d $)   iff any $x$ in $K(Z)$ can be written as the
sum of two products
\begin{equation}\label{eq101}
x = x_1 x_2\ldots x_d + y_1 y_2\ldots y_d
\end{equation}
with each $x_j$,  $y_j$ either in $K(A)$ or in $K(B)$,  with the first terms 
$x_1$  in $K(B)$ and $y_1$   in $K(A)$, and satisfying
$$   \prod\limits^d_1 \|x_j\|+\prod\limits^d_1 \|y_j\| \le C \|x\|_{  K(Z)
}.$$ 
 Let $D\ge 1$ be another integer and let $d=2D+1$. 
 We will say that
$L^A(Z)\le  D$ if there is a constant $C$ such that the same as before holds but 
with  
$x = x_1 x_2\ldots x_d$ or equivalently with the $y_j$'s all vanishing.
More precisely, we have $x_{2j+1}\in K(B)$ for $j=0,1,2...,D$ and  
$x_{2j}\in K(A)$ for $j=1,2...,D$.
For convenience, we say that
$L^B(Z)\le  D$ if there is a constant $C$ such that the same as before holds but 
with  
$x = y_1 x_2\ldots y_d$ or equivalently with the $x_j$'s all vanishing.
 
Finally, if $d=2D+1$ and if the property used above
 to define $L(Z)\le d$ holds for $n=1$ 
but with the $y_j$'s all vanishing, we say $L^A_1(Z)\le D$. Again we say
 $L^B_1(Z)\le D$ if this holds for $n=1$ but   with the $x_j$'s all vanishing.
  Needless to say $L(Z)=L(Z;A,B) $ is defined as the 
smallest integer $d>0$ such that $L(Z)\le d$, and similarly for 
$L^A(Z)$, $L^B(Z)$, $L_1^A(Z)$, and so on.

Roughly, $L(Z)\le d$ corresponds to factorizations of length at most $d$ jointly in
$K(A)$ and $K(B)$, while $L^A(Z)\le D$ corresponds to factorizations of  length
$D$ in $K(A)$, (and a fortiori of length at most $2D+1$ jointly in $K(A)$ and $K(B)$).

\begin{rem}\label{rem9}  By an elementary counting argument, we find:
$$  2 L^A(Z)-1\le L(Z) \le  2 L^A(Z)+1$$
$$   2L_1^A(Z)-1\le L_1(Z) \le  2 L_1^A(Z)+1.$$
Moreover (this is obvious if $A$ is   unital, otherwise   we   use an
approximate unit)
$$ L^A(Z)\le L^B(Z)+1  \quad{\rm and} \quad L_1^A(Z)\le L_1^B(Z)+1.$$
\end{rem}

\begin{rk} Note that length $\le d$ obviously passes to quotients: for any ideal  $I
\subset Z$, we have 
\begin{equation}\label{LQ}
L(Z/I)\le L(Z).
\end{equation}
\end{rk}
Let $A,B$ be two $C^*$-algebras. Let
$A\overset{\cdot}{\ast} B$ be their (non-unital)  $C^*$-algebraic free product.
This is obtained by completion of the algebraic free product with respect to the
maximal $C^*$-norm on it.
 Let $ V_d\subset A \overset{\cdot}{\ast} B$ be the 
subspace generated by elements of the form $x_1x_2\ldots x_d$ with $x_k\in A$ or 
$x_k\in B$ in such a way that $x_k$ and $x_{k+1}$ do not belong to the same 
subalgebra ($A$ or $B$).

Similarly, let $V^A_d$ (resp.\ $V^B_d$) be the closed span of elements of the 
form $x_1x_2\ldots x_d$ as above but such that $x_1\in A$ (resp.\ $x_1\in B$). 
Note that $V_d$ is obviously the closure of $V^A_d + V^B_d$.

Given   $C^*$-subalgebras $A,B$ as above we denote by
$$ \overset{\cdot}Q_{Z}\colon \  A\overset{\cdot}{\ast} B \to Z$$
the (surjective) $*$-homomorphism canonically extending the inclusions
$A\subset Z$ and $B\subset Z$.

\begin{thm}\label{L}
Let $A\subset Z$ and $B\subset Z$ be $C^*$-subalgebras
 as above. Let $d\ge 1$ be an integer.
The following assertions are equivalent:
\begin{itemize}
\item[(i)] $L(Z;A,B)\le d$
\item[(ii)] The restriction of the canonical quotient map 
$\overset{\cdot} Q_Z\colon \ A \dot * B \to Z$ to  
$V_{\le d}$ is a complete surjection.

Moreover,   (i) or (ii)  implies
\item[(iii)] Every $x$ in $K(Z)$ can be written as a product
\[x=x_1\cdots x_{d+1}\]
with $x_1,\cdots,x_{d+1}$ either in $K(A)$ or in $K(B)$.
\end{itemize}

\end{thm}
For completeness, we also state the analogue for $L_1$:

\begin{thm}\label{L1}
Let $A\subset Z$, $B\subset Z$ and $d\ge 1$ be as above.
The following assertions are equivalent:
\begin{itemize}
\item[(i)] $L_1(Z;A,B)\le d$
\item[(ii)] The restriction of the canonical quotient map 
$\overset{\cdot} Q_Z\colon \ A \dot * B \to Z$ to  
$V_{\le d}$ is a   surjection.

Moreover,   (i) or (ii)  implies
\item[(iii)] Every $x$ in $Z$ can be written as a product
\[x=x_1\cdots x_{d+1}\]
with $x_1,\cdots,x_{d+1}$ either in $K(A)$ or in $K(B)$,
with the understanding that
$x_1$ is a $1\times \infty$ matrix and $x_{d+1}$
is a $\infty\times 1$ matrix.
\end{itemize}

\end{thm}

To prove these statements, we will use the ``Haagerup tensor product"
of operator spaces, for which we refer the reader to \cite{P1,BLM,ER}.
The main relevant fact for our purpose is the following.
\begin{lem}\label{lemm3}

The space $V^A_d$ (resp.\ $V^B_d$) is completely isomorphic to the Haagerup 
tensor product 
$$A\otimes_h B \otimes_h A\ldots\ ({\rm resp.}\ B\otimes_h A \otimes_h 
B\ldots)$$
 with a total of $d$ factors. 
 \end{lem} 
 \begin{proof}
 For this last fact (apparently due to the author), we refer the reader  to 
\cite[Exercise 5.8, p.~108 and p.~433-434]{P1}. This is a refinement of results originally  in \cite{CES}.
\end{proof}

\begin{rem}\label{rem30}
It will be convenient to use the universal $C^*$-algebra generated by two 
projections $p,q$, denoted by $C_2$. We define $C_2$ as follows:\ let $x$ be a 
formal linear combination of the set
\[
J = \{1,p,q,(pq)^j, (qp)^j, (pq)^jp, (qp)^jq\mid j\ge 1\}.
\]
We set $\|x\|$ equal to the supremum of the norm of $x$ in $B(H)$ when we 
replace 
$p,q$ by an arbitrary pair of orthogonal projections in $B(H)$, $H$ being an 
arbitrary Hilbert space. Then $C_2$ can be defined as the completion of the 
space of these $x$'s equipped with this norm. Let $\vp_1 = p-(1-p) = 2p-1$ and 
$\vp_2 = 2q-1$. Note that $\vp_1,\vp_2$ are unitaries with  $\vp^2_1 = \vp^2_2 = 
1$, which generate $C_2$ as a $C^*$-algebra. Consequently (see \cite{RS}  for more on 
this), $C_2$ can be identified with $C^*({\bb Z}_2* {\bb Z}_2)$ the 
$C^*$-algebra of the (amenable) dihedral group, with $\vp_1$ and $\vp_2$  corresponding
to the  free generators of the two (free) copies of $  {\bb Z}_2 $. For convenience, we introduce 
the following notation:
\begin{alignat*}{2}
p_{2j} &= (pq)^j,&\quad p_{2j+1} &= (pq)^jp\\
q_{2j} &= (qp)^j,&\quad q_{2j+1} &= (qp)^jq.
\end{alignat*}
We will use the observation that the family $J$ is linearly independent in $C_2$. This 
is easy to check by observing that
\begin{align*}
(pq)^j &= (\vp_1\vp_2)^j + \text{lower order terms,}\\
(qp)^j &= (\vp_2\vp_1)^j + \text{lower order terms},
\end{align*}
and similarly for $(pq)^jp$ and $(qp)^jq$.

This observation implies that, if we fix $d\ge 1$, there 
is a constant $K(d)$ such that for any finitely supported families of scalars 
$(\lambda_j)_{j\ge 1}$ and $(\mu_j)_{j\ge 1}$ we have 
\begin{equation}\label{eq19}
K(d)^{-1} \max\{\sup_{j\le d}|\lambda_j|, \sup_{j\le d} |\mu_j|\} \le 
\left\|\sum_{j\le d} \lambda_jp_j + \sum_{j\le d} \mu_jq_j\right\| \le K(d) 
\max\{\sup_{j\le d} |\lambda_j|, \sup_{j\le d} |\mu_j|\}.
\end{equation}
\end{rem}

 Actually, the sum $V^A_d + V^B_d$ is a direct sum of operator spaces. More
precisely:

\begin{lem}\label{lem31}
Let $V_{\le d}$ be the closed span of $V_j$ for $1\le j\le d$. Then we have for 
each fixed $d\ge 1$ the following complete isomorphisms:
\begin{align}
\label{eq20}
V_{\le d} &\simeq V_1 \oplus\cdots\oplus V_d,\\
\label{eq21}
V_d &\simeq V^A_d \oplus V^B_d\\
\intertext{and consequently}
\label{eq22}
V_{\le d} &= V^A_1 \oplus V^B_1 \oplus\cdots\oplus V^A_d \oplus V^B_d.
\end{align}
\end{lem}

\begin{proof}
The proof is elementary. Consider an element $x$ in $V^A_1 + V^B_1 +\cdots+ 
V^A_d + V^B_d$, say we have $x = \alpha_1+\beta_1 +\cdots+ \alpha_d+\beta_d$ 
with $\alpha_j\in V^A_j$, $\beta_j \in V^B_j$. Let $p,q$ and $C_2$ be as in
the above  Remark \ref{rem30}. 
We may then consider the pair of (non-unital) representations $\pi_p\colon \ 
A\to (A\overset{\cdot}{\ast} B) \otimes_{\min} C_2$ and $\pi_q \colon \ 
B\to (A\overset{\cdot}{\ast} B) \otimes_{\min} C_2$ defined by $\pi_p(a) = 
a\otimes p$ and $\pi_q(b) = b\otimes q$. Let $\pi\colon \ A\overset{\cdot}{\ast} 
B 
\to (A\overset{\cdot}{\ast} B) \otimes_{\min} C_2$ be the representation 
canonically extending (jointly) $\pi_p$ and $\pi_q$. 

 Note that
\begin{equation}\label{eq222}
\pi(x) = \sum^d_{j=1} \alpha_j \otimes p_j + \beta_j\otimes q_j.
\end{equation}
By \eqref{eq19}, the span of $\{p_j,q_j\mid 1\le j\le d\}$ is (completely) 
isomorphic to ${\bb C}^{2d}$, therefore \eqref{eq22} follows immediately from \eqref{eq222}. Then 
\eqref{eq21} follows by restricting to $x=\alpha_d+\beta_d$ and \eqref{eq20} is 
but a combination of \eqref{eq21} and \eqref{eq22}.
\end{proof}

 \begin{proof}[Proof of Theorems \ref{L}  and \ref{L1}]
To simplify the notation, let us denote  
by $X_d=A\otimes_h B \otimes_h A\otimes_h \ldots$ (resp. $Y_d=B\otimes_h A \otimes_h 
B\ldots$) 
where each   $X_d$ and $Y_d$ have exactly $d$ factors.
By definition of the Haagerup tensor product  the assumption
that $L(Z)\le d$ (resp.\ $L^A(Z)\le d$) equivalently means that the product map
$X_d+Y_d \to Z$ (resp. the product map $Y_{2d+1} \to Z$)  is a
complete surjection (see e.g.
\cite[Cor. 5.3 p. 91]{P1}). Here $X_d+Y_d$ is defined as the operator  quotient space
of the direct sum (say in the $\ell_1$-sense) $X_d\oplus Y_d$ by the kernel of the
map $(x,y)\to x+y$ (see \cite[p. 55]{P1} for more information). Since we just saw 
(by \eqref{eq21} and Lemma \ref{lemm3}) that $V_d$  (resp. $V^B_{2d+1}$) can be
identified with $X_d \oplus Y_d$ (resp. with $Y_{2d+1}$), this holds iff   the
restriction of $\overset{\cdot}Q_{Z}$  to $V_d$ (resp.\ $V^B_{2d+1}$)  is a
complete surjection.
 The same argument yields the analogous statement concerning
$L_1(Z)\le d$ (resp.\ $L_1^A(Z)\le d$). Finally, the assertions (iii)
are proved using a unit if it exists, an approximate unit otherwise.
 \end{proof}

\begin{rk} Similarly, $L^B(Z)\le d$ iff $\overset{\cdot}Q_{Z}$ restricted
to 
$V^A_{2d+1}$ is a complete surjection.
\end{rk}

 The analogous notation and statements for the {\it unital} free product are as follows.
 
\noindent Let 
$A {\ast} B$ be the  unital free product of $A$ and $B$ (both assumed unital).
Clearly, there is a canonical surjective $*$-homomorphism
$\kappa\colon \ A\overset{\cdot}{\ast} B \to A {\ast} B$.
 Let 
$Q_Z\colon \ A {\ast} B \to Z$ be the natural 
(quotient) unital $*$-homomorphism. We have obviously $Q_{Z}\
\kappa=\overset{\cdot}Q_{Z}$.

Let $W_{\le d}$ be the subspace generated by elements of the form $x_1x_2\ldots 
x_d$ with either $x_k$ in $A$ or $x_k$ in $B$ for each each $k=1,\ldots, d$.  
Let $\varphi$ (resp.\ $\psi$) be a state on $A$ (resp.\ $B$).
 Let 
$\overset{\circ}{A}$ (resp.\ $\overset{\circ}{B}$) denote the subspace of $A$ 
(resp.\ $B$) formed of all elements on which $\varphi$ (resp.\ $\psi$) vanishes. 

We will   keep this choice of states $\varphi$ (resp.\ $\psi$) fixed throughout
the rest of the paper. Note that   
$\overset{\circ}{A}$ (resp.\ $\overset{\circ}{B}$) implicitly depend
on this initial choice, even though the notation does not indicate it.

Note that $A\simeq {\bb C} 1_A \oplus \overset{\circ}{A}$ and $B\simeq {\bb C} 
1_B \oplus \overset{\circ}{B}$.  We denote by $W^A_d$ (resp.\ $W^B_d$) the 
closed span in $A\ast B$ of all elements $y$ of the form
\begin{equation}\label{eq33}
y = y_1y_2\ldots y_d
\end{equation}
with each $y_k$ either in $\overset{\circ}{A}$ or in $\overset{\circ}{B}$ in such a 
way that no two consecutive elements   belong to the same set $\overset{\circ}{A}$ 
or $\overset{\circ}{B}$ (so if $y_k\in \overset{\circ}{A}$ then $y_{k+1}\in 
\overset{\circ}{B}$) and finally such that $y_1\in \overset{\circ}{A}$ (resp.\ 
$y_1\in \overset{\circ}{B}$).

Roughly speaking, $W^A_d$ (resp.\ $W^B_d$) is spanned by the elements of length 
exactly equal to $d$, that start in $A$ (resp\ $B$). We denote by $W_d$ the 
closure of $W^A_d + W^B_d$.

Let us denote by ${\cl W}^A_d$ (resp.\ ${\cl W}^B_d$) the linear span in $A\ast B$ 
of all elements $y$ of the form \eqref{eq33} with $y_1$ in $A$ (resp.\ $y_1$ in 
$B$).
Let $${\cl W}_d = {\cl W}^A_d + {\cl W}^B_d.$$
 With this notation, $W^A_d$ (resp.\ $W^B_d$) appears as the closure of 
${\cl W}^A_d$ (resp.\ ${\cl W}^B_d$) in $A{\ast} B$, 
and $W_d$ is the closure of ${\cl W}_d$.
 
Note that ${\cl W}^A_d$ (resp.\ ${\cl W}^B_d$) is clearly linearly isomorphic to 
the algebraic tensor product $\overset{\circ}{A} \otimes \overset{\circ}{B} 
\otimes \cdots$ (resp.\ $\overset{\circ}{B} \otimes \overset{\circ}{A}\otimes 
\cdots$). 
Therefore, using the canonical embeddings $A\subset A\overset{\cdot}{\ast} B$ and $B\subset A 
\overset{\cdot}{\ast}B$ we can unambiguously define linear embeddings of ${\cl W}^A_d$ and ${\cl 
W}^B_d$ into $A\overset{\cdot}{\ast}B$.
This gives us a linear embedding of ${\cl W}_d$ into $A\overset{\cdot}{\ast}B$. 
Let us denote by 
$\Lambda$ the linear map extending the preceding one to the linear span of 
$\{{\cl W}_j\mid j\ge 1\}$. Then we have
 
\begin{lem}\label{lem32}
For any $d\ge 1$, the mapping $\Lambda$ extends to a complete isometry from
$W_1+...+W_{  d}$ into $A\overset{\cdot}{\ast} B$, which lifts the canonical
surjective $*$-homomorphism 
$\kappa\colon \ A\overset{\cdot}{\ast}B \to A{\ast} B$.
\end{lem}

\begin{proof}
Let $x = \omega_1 +\cdots+ \omega_d$ with $\omega_j \in {\cl W}_j$ for any $j\ge 
1$. We will show that $\|{\Lambda}(x)\|_{A\dot *B} = \|x\|_{A{\ast} B}$. Note that 
${\Lambda}$ is  
trivially a lifting of $\kappa$, so that  
$
q({\Lambda}(x)) = x$ and hence  $\|x\|\le \|{\Lambda}(x)\|
$
is immediate. To prove the converse, consider a pair of representations 
\[
\pi_1\colon \ A\to B(H)\quad \text{and}\quad \pi_2\colon \ B\to B(H)
\]
such that the   associated representation $\pi$    on $A\overset{\cdot}{\ast}
B$ is isometric, so that
$\|\pi({\Lambda}(x))\| = \|{\Lambda}(x)\|$. Let $\pi_1(1_A) = p$ and $\pi_2(1_B)=q$. We may 
introduce the maps $\widehat\pi_1$ and $\widehat\pi_2$ on $A$ and $B$ respectively by 
setting
\begin{align*} 
\widehat\pi_1(a) &= \pi_2(a) + \varphi(a)(1-p)\\
\widehat\pi_2(b) &= \pi_2(b) + \psi(b)(1-q).
\end{align*}
Note that since $\pi_1(a) = p\pi_1(a)p$ and $\pi_2(b) = q\pi_2(b)q$, $\widehat\pi_1$ 
and $\widehat\pi_2$ are unital completely positive (c.p.\ in short) maps. By 
\cite{Bo1} (see also \cite{Bo2}), there is a unital completely positive map 
$\widehat\pi\colon \ A{\ast} B\to B(H)$ such that for any $y$ as in \eqref{eq33} we have 
\begin{align*}
\widehat\pi(y) &= \widehat\pi_1(y_1) \widehat\pi_2(y_2)\ldots \text{ if } y_1 \in 
\overset{\circ}{A}\\
\intertext{and}
\widehat\pi(y) &= \widehat\pi_2(y_1) \widehat\pi_1(y_2)\ldots \text{ if } y_1\in 
\overset{\circ}{B}.
\end{align*}
But for any $y_1$ in $\overset{\circ}{A}$ (resp.\ $y_2$ in $\overset{\circ}{B}$) 
we have $\widehat\pi_1(y_1) = \pi_1(y_1)$ (resp.\ $\widehat\pi_2(y_2) = \pi_2(y_2)$).

Hence this shows that $\widehat\pi(\omega_j) = \pi({\Lambda}(\omega_j))$ for any $j$ and 
hence
\[
\widehat\pi(x) = \pi({\Lambda}(x)).
\]
Thus we conclude that
\[
\|\pi({\Lambda}(x))\| \le \|\widehat\pi\| \|x\| \le \|x\|
\]
and since we choose $\pi$ so that $\|\pi({\Lambda}(x))\| = \|{\Lambda}(x)\|$, we obtain as 
announced $\|{\Lambda}(x)\| \le \|x\|$.
This shows that ${\Lambda}$ is isometric. The proof that it is completely isometric is 
entirely similar, we leave the routine details to the reader.
\end{proof}
 
We have then

\begin{lem}\label{lem33}
For each fixed $d\ge 1$, we have the following complete isomorphisms:
\begin{align}
\label{eq23}
W_{\le d} &\simeq {\bb C} \oplus W_1 \oplus\cdots\oplus W_d\\
\label{eq24}
W_d &\simeq W^A_d \oplus W^B_d\\
\label{eq25}
W_{\le d} &\simeq {\bb C} \oplus W^A_1 \oplus W^B_1 \oplus\cdots\oplus W^A_d 
\oplus W^B_d.
\end{align}
Moreover, we have complete isomorphisms
\[
W^A_d \simeq \overset{\circ}{A} \otimes_h \overset{\circ}{B} \otimes_h \cdots 
\quad \text{and}\quad W^B_d \simeq \overset{\circ}{B} \otimes_h
\overset{\circ}{A} \otimes_h\cdots~.
\]
\end{lem}
\begin{proof} 
Let us first check that the sums $W_1 \oplus\cdots\oplus W_d$  and   $W^A_d \oplus W^B_d$ are direct ones. Using the lifting in Lemma \ref{lem32}, this is an immediate consequence of
\eqref{eq20} and \eqref{eq21}, since we have,  of course, ${\Lambda}(W_d)\subset V_d$, 
${\Lambda}(W_d^A)\subset V_d^A$ and ${\Lambda}(W_d^B)\subset V_d^B$. In particular, this   proves  \eqref{eq24}.
Let $\phi\ast \psi$ denote the free product state on $A\ast B$ (see \cite[p. 4]{VDN}).
Note that  $\phi\ast \psi$ vanishes on $W_1+...+W_d$ and hence \eqref{eq23} follows.  
 Then \eqref{eq25} 
is but a recapitulation. Finally, the last assertion follows
from the Lemmas \ref{lemm3} and  \ref{lem32} using again  ${\Lambda}(W_d^A)\subset V_d^A $
and ${\Lambda}(W_d^B)\subset V_d^B $, and the injectivity of the Haagerup tensor product (cf. e.g. \cite[p. 93]{P1}).
\end{proof}
 
In the next statement, we denote by $W^B_{\le 2d+1}$  
 the sum $\sum_{j=1}^{2d+1}W^B_{j}$. Note that the latter sum is closed
 since, by Lemma \ref{lem33} it is a direct sum. Equivalently, this is the closed span
 of all  alternated products in $\overset{\circ}{A}$ and $\overset{\circ}{B}$ but
with at most $d$ factors in $\overset{\circ}{A}$.

\begin{thm}\label{LU}
Let $A,B,Z$ be unital $C^*$-algebras as above. Then, $L(Z)\le
d$  (resp.\ $L^A(Z)\le d$) iff the restriction of $Q_{Z}\colon \ 
A{\ast} B \to Z$ to $W_{\le d}$ (resp.\ $W^B_{\le 2d+1}$) 
is a complete surjection.
 Moreover, $L_1(Z  )\le d$ (resp.\ $\ell^A_1(Z)\le 
d$) iff the restriction of $Q_{Z}$ to $W_{\le d}$ (resp.\ $W^B_{\le 2d+1}$)   is a
surjection.
\end{thm}
\begin{proof}
Assume  $L(Z)\le d$. Then since  $Q_{Z}\ \kappa=\overset{\cdot}Q_{Z}$, 
Theorem 
 \ref{L} implies that $Q_{Z}$ restricted to $\kappa(V_d)$ is a complete
surjection. Since
$\kappa(V_d)\subset W_{\le d}$, it follows that $Q_{Z}$ restricted to  $W_{\le d}$ is
a complete surjection. Conversely, assume that $Q_{Z}$ restricted to  $W_{\le d}$
is a complete surjection. Since $W_{\le d}$ is spanned by the unit of $A{\ast} B   $
and $W_1+...+W_d$, Lemma
\ref{lem32}, recalling  \eqref{eq23}, implies that there is a completely bounded
map $\hat {\Lambda}\colon\ W_{\le d}\to A\overset{\cdot}{\ast}  B$ lifting $\kappa$,
 defined e.g. by
$\hat {\Lambda}(\lambda 1 + x)=\lambda 1_A+{\Lambda}(x)$.
Note that $\hat {\Lambda}(W_{\le d})\subset V_{\le d}$ and hence 
$\overset{\cdot}Q_{Z}$
 restricted to  $V_{\le d}$ is a complete surjection.
 
 Now observe that  the completely isometric   mapping
 $$x\to 1_A x \quad ({\rm resp.}\quad x\to 1_Bx)$$
 takes $V_{d-1}^B$ (resp.  $V_{d-1}^A$) to  $V_{d}^A$  (resp.  $V_{d}^B$)
 and we have
 $Q_{Z}(x)= Q_{Z}(1_A x)= Q_{Z}(1_B x)$,  and similarly for $V_{d-2}^B$ (resp. 
$V_{d-2}^A$),
 $V_{d-3}^B$ (resp.  $V_{d-3}^A$), and so on. 
Since, by Lemma \ref{lem31},  $V_{\le d}$
 decomposes as a direct sum of $V_{j}^A+V_{j}^B$ ($j\le d$), 
it is easy to use the preceding observation
 to replace the elements of $V_{\le d}$ by suitably chosen ones 
in $V_{ d}$ in order to show
 that $\overset{\cdot}Q_{Z}$ restricted to  the smaller subspace $V_d\subset V_{\le
d}$ is a complete surjection.
 By Theorem \ref{L} again, we conclude that  $L(Z)\le d$.
 This proves the part of the statement concerning $L(Z)$, but actually
 the same proof also establishes the   part concerning $L_1(Z)$ by removing
``complete"
 from ``complete surjection". The other part is proved similarly.
 We leave the details to the reader. 
 \end{proof}

 Finally, we give the basic result that connects the length
with the first part of the paper.

\begin{lem}\label{lem22} Assume $L(Z)\le d$. Then there is a constant
$C$ such that for any
bounded homomorphism $\Phi\colon\ Z\to B(H)$, we 
have
\[\|\Phi\|_{cb} \le C(\max\{ \|\Phi_{|A}\|_{cb},  \|\Phi_{|B}\|_{cb} \})^d .
\]
Moreover for any bounded derivation $\Delta\colon\ Z\to B(H)$
(relative to a representation of $Z$ on $B(H)$)
\[\|\Delta\|_{cb} \le C(d+1)\max\{ \|\Delta_{|A}\|_{cb},  \|\Delta_{|B}\|_{cb} \})  .
\]
\end{lem}
\begin{proof}
 Let $t=\max\{ \|\Phi_{|A}\|_{cb},  \|\Phi_{|B}\|_{cb} \}$.
Let $\Phi_n=Id\otimes \Phi \colon\ M_n(Z) \to M_n(B(H))$.
Let $x\in M_n(Z)$ with $\|x\|_{  M_n(Z)}<1$. With the notation in \eqref{eq100} we
have 
\[ \| \Phi_n(x-x_1 x_2\ldots x_d - y_1 y_2\ldots y_d   )  \| \le \|  \Phi_n \| \vp.\]
But we have clearly
\[ \| \Phi_n( x_1 x_2\ldots x_d + y_1 y_2\ldots y_d   )  \| \le t^d (\prod\limits^d_1
\|x_j\|+\prod\limits^d_1 \|y_j\|) \le Ct^d  \]
hence 
$\| \Phi_n(x     )  \| \le Ct^d +\|  \Phi_n \| \vp$ and hence
letting $\vp \to 0$ (here we crucially use that $\Phi$ is continuous!)   
  $\| \Phi_n\|\le Ct^d$ and  $\|\Phi\|_{cb} \le Ct^d$.
A similar argument gives the other inequality.
\end{proof}

The main result of \cite{P4}  says that essentially we have a converse
(note however that this is not exactly the converse, see Remark \ref{insist} below).

\begin{thm} 
Let $A\subset Z, B\subset Z$ as before. Assume $A,B,Z$ unital with unital embeddings.
Let $Alg(A,B)$ denote the (dense) subalgebra generated by
$A$ and $B$. If there is a constant $C$ such that any homomorphism
$\Phi\colon\ Alg(A,B)\to B(H)$ satisfies 
$\|\Phi\|_{cb} \le C(\max\{ \|\Phi_{|A}\|_{cb},  \|\Phi_{|B}\|_{cb} \})^d ,
$
then $$L(Z)\le d.$$

\end{thm}
\begin{proof} This is a particular case of Theorem 6 in \cite{P4}. 
\end{proof}

\begin{rk} For simplicity, in the definition of length and in  Theorems \ref{L}  and
\ref{LU}, we have restricted our description to pairs $A,B$ of subalgebras, but it is
easy to extend these statements (by a simple iteration) to triples $A,B,C$ of
subalgebras, or to any finite given number $N$ of them. Of course, the resulting
(complete) isomorphism constants in Lemmas \ref{lem31}  and \ref{lem33} will now
depend both on
$d$ and $N$.
\end{rk} 
\begin{rk} Throughout this section, we have restricted attention to
$C^*$-algebras but it is easy to verify that the same results remain valid
when $A,B,Z$ are non self-adjoint operator algebras with minor changes in the
proofs.   

\end{rk} 
\section*{Length for the maximal tensor product}

In this section, we will specialize the 
preceding  to the situation when
$Z=A\otimes_{\max} B$, where $A,B$ are two unital $C^*$-algebras
embedded into $Z$ via the mappings
$a\to a\otimes 1$ and $b\to 1\otimes b$. We will identify
$A$ with $A\otimes 1$ and $B$ with $1\otimes B$, and view them as subalgebras
of $Z=A\otimes_{\max} B$.

We will denote for simplicity
$$L(A\otimes_{\max} B)=L(A\otimes_{\max} B; A\otimes 1, 1\otimes B)$$
and similarly for 
$L_1(A\otimes_{\max} B), L^A_1(A\otimes_{\max} B) $ and $L^B_1(A\otimes_{\max} B) $.

We will similarly denote
$L(A\otimes_{\min} B)=L(A\otimes_{\min} B; A\otimes 1, 1\otimes B)$, and
similarly for $L_1,L^A_1,L^B_1$.

In \cite{P4}, the author introduced the ``similarity degree"  a pair of unital 
 $C^*$-algebras $A,B$ 
as follows:\ 

\noindent Changing the notation from  \cite{P4} slightly, we will say here
that $d(A,B)\le d$
(resp.   $d_1(A,B)\le d$) if there is a constant $C$ such that 
for any pair $u\colon \ A\to B(H)$, $v\colon \ B\to B(H)$ of c.b.\ unital homomorphisms 
with commuting ranges, the homomorphism $u\cdot v\colon\ A\otimes B\to B(H)$ 
(taking $a\otimes b$ to $u(a)v(b)$) is c.b.\ 
(resp. is bounded) on $A\otimes_{\max} B$ with c.b.\ 
norm (resp. with norm) $\le C \max(\|u\|_{cb}, \|v\|_{cb})^d$.

Of course, the number $d(A,B)$ (resp.   $d_1(A,B)$)  is defined as the infimum of the numbers $d\ge 1$ such that this property  holds.

If $(A,B)$ is such that for any $(u,v)$ as above, the map $u\cdot v$ is c.b.\ (resp. bounded) on 
$A\otimes_{\max} B$, then $d(A,B)<\infty$ 
 (resp.   $d_1(A,B)<\infty$) (see \cite{P5}).
 
 Let us denote by $d^ A(A, B)$ 
(resp. $d_1^ A(A, B)$) the smallest $d$ with the following property:
there is a constant $C$ such that
for any c.b. unital homomorphism $u:\ A \to B(H)$ and any unital $*$-homomorphism $\sigma:\ B\to B(H)$
with commuting ranges (i.e. we have   $\sigma(B)\subset u(A)'$), the product
mapping $u.\sigma$ is c.b. 
(resp. bounded) on $A\otimes_{\max} B$ with c.b. norm  (resp. with norm)
$\le C\|u\|^d_{cb}$.  
Moreover, we set by convention
$$d^ B(A, B)   = d^ B(B, A)\quad d_1^ B(A, B)   = d_1^ B(B, A) .$$
When $A$ is nuclear, we claim that this holds with $K=1$ and $d=2$ and then  it even  
holds for all {\it complete contractions} $\sigma:\ B\to u(A)'$. Indeed, if $A$ is nuclear, 
for any operator space $B$,  the
mapping $q:\ A\otimes_h B \otimes_h A\to A\otimes_{\min} B$ 
defined by  $q(a\otimes b\otimes a')= aa'\otimes b$ is a complete metric surjection (see \cite[p. 240-241]{P2}) .
Let $P_3 :\ B(H) \otimes_h B(H) \otimes_h B(H)\to B(H)$ 
be the product map, which is clearly a complete contraction.
We have obviously
$$u.\sigma (aa'\otimes b)=(u.\sigma )q (a\otimes b\otimes a')=P_3 (u\otimes \sigma\otimes u)(a\otimes b\otimes a').$$
Hence since $q$ is a  complete metric surjection, we have
$$\|u.\sigma\|_{CB(  A\otimes_{\min} B , B(H)) }= \|(u.\sigma)q\|_{cb}=  \|P_3 (u\otimes \sigma\otimes u)\|_{cb} \le \|u\|_{cb}^2 \|\sigma\|_{cb} .$$
This holds for any operator space $B$.
A fortiori, when $B$ is a $C^*$-algebra, 
we can replace the min-norm by the max-one and this proves the above  claim.

 \begin{rk} In \cite{P5}, we introduced the similarity degree 
$d(A)$ of a $C^*$-algebra $A$. This is defined
 as the smallest $d\ge1$ such that there is a constant $C$ so that 
  any bounded homomorphism
 $u\colon\ A\to B(H)$ satisfies $$\|u\|_{cb}\le C\|u\|^d.$$ This is related to the number $d(A,B)$ via
 the following obvious estimate:
 $$d( A\otimes_{\max} B  )  \le d(A,B) \max\{d(A),d(B)\}.$$
 
 \end{rk}

 As explained in \cite{P4}, the number $d(A,B)$ and its other
 variants are closely related to the length 
  $L(A \otimes_{\max} B)$. 
Let us briefly recall this here.  

For a pair $(A,B)$ of $C^*$-algebras, we will consider the following 
properties (see also \cite{LM}):
\begin{itemize}
\item[]
\begin{itemize}
\item[(SP)] For any pair $u\colon \ A\to B(H)$, $v\colon \ B\to B(H)$ of c.b.\ 
homomorphisms with commuting ranges, the product map $u\cdot v$ is c.b.\ 
on $A\otimes_{\max}B$.
\item[(SP)$_1$] For any pair $(u,v)$ as in (SP), the product map $u\cdot v$ is 
bounded on $A\otimes_{\max} B$.
\end{itemize}
\end{itemize}
Then the main result concerning $d(A,B)$ in \cite{P4} can be stated as follows.

\begin{thm}[\cite{P4}]
Assume (SP) (resp.\ (SP)$_1$). Then necessarily 
$d(A,B)<\infty$ (resp.\ $d_1(A,B)<\infty$), and moreover 
\[
d(A,B) = L(A\otimes _{\max} B)\quad (\text{resp. } d_1(A,B) = L_1(A\otimes _{\max} B)).
\]
\[
d^A(A,B) = L^A(A\otimes _{\max} B)\quad (\text{resp. } d_1^A(A,B) = L_1^A(A\otimes _{\max} B)).
\]
\end{thm}

\begin{rem}\label{insist} Assuming (SP)$_1$, we get by Lemma \ref{lem22}   that
$d(A,B)\le L(A\otimes _{\max} B)$ and also $d_1(A,B)\le L_1(A\otimes _{\max} B)$
or $d^A(A,B)\le L^A(A\otimes _{\max} B)$ and 
$d_1^A(A,B)\le L_1^A(A\otimes _{\max} B)$. However,
it may be worthwhile to insist on an unpleasant feature of this particular 
setting involving $A\otimes_{\max} B$: If we only assume
 $L(A\otimes _{\max} B)<\infty$ (resp.\ 
$L_1(A\otimes _{\max} B)<\infty$), we cannot verify in full generality that
\[
d(A,B)\le L(A\otimes _{\max} B)\quad (\text{resp. } d_1(A,B)\le L_1(A\otimes _{\max} B))
\]
because we do not see how to check that (SP) (resp. (SP)$_1$) holds. The difficulty
lies in the  fact that we have an approximate factorization in \eqref{eq100} relative
to a  norm (the max-norm) for which we do not know yet that $u\cdot v$ is continuous
!  See the proof
of  Lemma \ref{lem22} for clarification. An equivalent difficulty arises with
\eqref{eq101}.  Fortunately, in  the situations of interest to us,  (SP)$_1$ holds
(or $u\cdot v$ is continuous) so there is no problem.
\end{rem}

\begin{rk} Note that (iii) in Theorem \ref{thm1} means that that $(A,B)$ satisfies
(SP) for all $B$. Note also that this is formally equivalent to saying that $(A,B)$
satisfies (SP)$_1$ for all $B$. Indeed, the latter
implies the existence of a function $F$
such that $\|u\cdot v\colon A\otimes _{\max} B \to B(H)\|\le
F(\|u\|_{cb},\|v\|_{cb})$ holds for all $B,u,v$. We can then show
that $\|u\cdot v\colon A\otimes _{\max} B \to B(H)\|_{cb} \le
F(\|u\|_{cb},\|v\|_{cb})$     by replacing $B$ by
$M_n(B)$ to estimate the cb-norm of $u\cdot v$. Thus $A$ is nuclear iff
for any $B$ the pair $(A,B)$ satisfies (SP)$_1$.

\end{rk}
 
Given unital $C^*$-algebras $A$ and $B$ we denote by
$$ \overset{\cdot}Q_{A,B}\colon \  A\overset{\cdot}{\ast} B \to A \otimes_{\max} B$$
the (surjective) $*$-homomorphism canonically extending the inclusions
$A\subset A \otimes_{\max} B$ and $B\subset A \otimes_{\max} B$.
 
The next two statements recapitulate what we know
from Theorems \ref{L} and \ref{LU}.

\begin{pro}\label{pro14}
Let $A,B$ be unital $C^*$-algebras satisfying (SP). Then, $L(A\otimes _{\max} B)\le d$ 
(resp.\ $L^A(A\otimes _{\max} B)\le d$) iff the restriction of $\overset{\cdot}Q_{A,B}\colon \ 
A\overset{\cdot}{\ast} B \to A \otimes_{\max} B$ to $V_d$ (resp.\ $V^B_{2d+1}$) 
is a complete surjection. Moreover, $L_1(A\otimes _{\max} B)\le d$ (resp.\
$L^A_1(A,B)\le  d$) iff the restriction of $\overset{\cdot}Q_{A,B}$ to $V_d$ (resp.\
$V^B_{2d+1}$)   is a surjection.
\end{pro}

\begin{pro}
Let $A,B$ be unital $C^*$-algebras satisfying (SP). Then, $L(A\otimes _{\max} B)\le d$ 
(resp.\ $L^A(A\otimes _{\max} B)\le d$) iff the restriction of $Q_{A,B}\colon \ 
A{\ast} B \to A \otimes_{\max} B$ to $W_{\le d}$ (resp.\ $W^B_{\le 2d+1}$) 
is a complete surjection. Moreover, $L_1(A\otimes _{\max} B)\le d$ (resp.\
$L^A_1(A,B)\le  d$) iff the restriction of $Q_{A,B}$ to $W_{\le d}$ (resp.\ $W^B_{\le
2d+1}$)   is a surjection.
\end{pro}

The preceding Theorem \ref{thm1} shows that (SP) holds for all $B$ iff $A$ is nuclear (and 
in that case $d(A,B)$ as defined above is $\le 3$, and also
 $d^A(A,B)\le 2$ and  $d^B(A,B)\le 1$).
Now we have  conversely:

\begin{thm}\label{thm6}
Let $A,B$ be $C^*$-algebras. If   $B$ is liberal   
and if $L(A\otimes _{\max} B)<\infty$
or if more generally $L(A\otimes _{\min} B)<\infty$,
 then $A$ is nuclear.
\end{thm}
 \begin{proof} Note $L(A\otimes _{\min} B)\le L(A\otimes _{\max} B)$.
Using Lemma \ref{lem22} with $Z=A\otimes _{\min} B$
and  recalling \eqref{eqH}, we see that $L(A\otimes _{\min} B)<\infty$
implies   property (vi) in Theorem \ref{thm1B}.
\end{proof}
As a corollary, we can answer several questions raised in \cite{P4}:

\begin{cor}\label{cor7}
Assume $\dim H=\infty$. Let $G$ be a discrete group. Then $L(C^*(G)\otimes _{\max}
B(H)) < 
\infty$ iff $G$ is amenable. Therefore, $L(C^*(\F_n)\otimes _{\max} B(H))=\infty$
for any $n\ge 2$.
 Moreover, $L(B(H)\otimes _{\max} B(H)) = \infty$.
\end{cor}

\begin{proof}
Recall that $C^*(G)  $ is nuclear iff $G$ is amenable (see \cite{L}).
 Moreover, by \cite{W}, $B(H)$ is
not nuclear.
Obviously $C^*(\F_\infty)$ is liberal. 
 A combination of \cite[p. 205  ]{A} and \cite{W} shows that $B(\ell_2)$ (or
 even  the Calkin algebra $B(\ell_2)/K (\ell_2)  $)
 is liberal. Therefore, this corollary follows from the preceding  theorem. 
\end{proof}

 As mentioned already in \cite{P4}, if $A$ is nuclear and $B$ an arbitrary $C^*$-algebra, then
$   L^A(A\otimes _{\max}B)\le 2$ and $   L(A\otimes _{\max}B)\le 3$.
Note the obvious inequality  $L(A\otimes _{\min}B)\le L(A\otimes _{\max}B)$.
 Then here is a final
  recapitulation:
\begin{thm}\label{thm8} Let $A$ be a $C^*$-algebra. The following are equivalent:
\begin{itemize}
\item[(i)] $A$ is nuclear.
\item[(ii)] For any $C^*$-algebra $B$, we have  $L(A\otimes _{\max}B)<\infty$.
\item[(ii)$_1$] For any $C^*$-algebra $B$, we have  $L_1(A\otimes _{\max}
B)<\infty$.
\item[(iii)]  For any $C^*$-algebra $B$, we have $ L^A(A\otimes _{\max}B)\le 2$.
\item[(iv)]  For any $C^*$-algebra $B$, we have $ L^B(A\otimes _{\max} B)\le 1$.
\end{itemize}
Moreover, these are all equivalent to the same properties with respect to
$A\otimes _{\min}B$.
\end{thm}

\begin{proof} The fact that nuclear implies (iii) or (iv)  
 (and a fortiori any of the other properties)
follows from properties of the so-called $\delta$-norm in \cite[p. 240]{P1}.
 The converses all follow
from the preceding theorem, recalling Remarks \ref{norm} and \ref{rem9}.
\end{proof}

\begin{rk} As already mentioned in Remark \ref{rem1}, 
it may be true that $L(A\otimes _{\max} B)<\infty$, 
or even merely $L_1(A\otimes _{\min} B)<\infty$, for a fixed pair $A,B$
implies $A\otimes_{\max} B= A\otimes_{\min} B$. One simple minded  approach to prove
this would be as follows: Let $Z= A\otimes_{\min} B$ and $n=1$ in \eqref{eq100}.
Consider $x\in A\otimes B$ ({\it
algebraic} tensor product) with $\|x\|_{\min}<1$. The problem
is simply to prove that the property
expressed by \eqref{eq100} (i.e. the fact that $L_1(A\otimes _{\min} B)\le d$)
automatically implies another representation as in \eqref{eq100} but for $\vp=0$.
Indeed, it is clear (here recall $n=1$) that  $\| \prod^d_1 x_j+\prod^d_1 y_j 
\|_{\max}<C$ so we would conclude
$\| \cdot  \|_{\max}\le C\| \cdot  \|_{\min}$. Note that if such a simple minded 
 proof (in particular not using \cite{Co})
is found, it would give a more direct way to
show   that nuclear passes to quotients.

\end{rk}

\noindent{\it Acknowledgement.} I am grateful to Ken Dykema for  a stimulating
question.

\end{document}